\documentclass{gtart}
\usepackage{rlepsf}
\input gtmonout
\volumenumber{1}
\volumeyear{1998}
\volumename{The Epstein birthday schrift}

\pagenumbers{295}{301}
\papernumber{13}
\received{15 November 1997}
\published{27 October 1998}

\def\<{\langle}
\def\>{\rangle}

\newtheorem{Th}{Theorem}[section]

\begin{document}

\title{On the fixed-point set of automorphisms of\\
non-orientable surfaces without boundary}

\shorttitle{Fixed-point sets of automorphisms of surfaces}

\author{M Izquierdo\\D Singerman}

\address{Department of Maths, M\"alardalen University,
721 23 V\"as\-te\-r\aa s, Sweden\\{\rm and}\\ 
Department of Maths,
University  of Southampton,
Southampton, SO17 1BJ, UK}

\asciiaddress{Department of Mathematics\\
Malardalen University\\
721 23 Vasteras, Sweden\\ 
Department of Mathematics\\
University  of Southampton\
Southampton SO17,1BJ UK}

\email{mio@mdh.se, ds@maths.soton.ac.uk}

\begin{abstract}
 Macbeath gave a formula for the number of fixed points for
 each non-identity element of a cyclic group of automorphisms of a compact
Riemann surface
in terms of the universal covering transformation group of the cyclic
group. We observe that
this formula generalizes to determine the fixed-point set of each
non-identity element of a cyclic
group of automorphisms acting on a closed non-orientable surface with one
exception; namely, when
this element has order 2. In this case the fixed-point set may have simple
closed curves
(called {\it ovals}) as well as fixed points. In this note we extend
Macbeath's results to include the
number of ovals and also determine whether they are twisted or not.
\end{abstract}

\asciiabstract{Macbeath gave a formula for the number of fixed points for
 each non-identity element of a cyclic group of automorphisms of a compact
Riemann surface
in terms of the universal covering transformation group of the cyclic
group. We observe that
this formula generalizes to determine the fixed-point set of each
non-identity element of a cyclic
group of automorphisms acting on a closed non-orientable surface with one
exception; namely, when
this element has order 2. In this case the fixed-point set may have simple
closed curves
(called ovals) as well as fixed points. In this note we extend
Macbeath's results to include the
number of ovals and also determine whether they are twisted or not.}




\maketitle

\cl{\small\it For David Epstein on the occasion of his sixtieth birthday}

\section{Introduction}

Let $Y$ be a compact non-orientable Klein surface of genus $p \ge 3$. By
genus here
we mean the number of cross-caps of the surface. Let $t\co  Y \to Y$ be an
automorphism
of order $M$. If $1\le i < M$ and if $i \neq M/2$ then the fixed-point set of
 $t^{i}$ consists of isolated fixed points and their number can be
calculated, as
 described below, by a formula which is completely analogous to Macbeath's
formula
 \cite{[M2]} concerning automorphisms of Riemann surfaces. However, if
$M=2N$ then the
fixed-point set of the involution $t^{N}$ consists of a finite number $n$ of
 disjoint simple closed curves called {\it ovals} together with a finite
number of
 isolated fixed points \cite{[BCNS]}, \cite{[Sch]}. The ovals may be {\it
twisted}
 or {\it untwisted} which means that they have M{\"o}bius band or annular
neigbourhoods
respectively.

In this note we calculate the number of ovals and isolated fixed-points of
$t^{N}$ and whether the ovals are twisted or not.

The information is given, as in Macbeath \cite{[M2]} in terms of the universal
 covering transformation group.

The authors acknowledge M{\"a}lardalen University and the Swedish Natural
Science Research Council for financial support.

\section{ The universal covering transformation group}

If $Y$ is a compact non-orientable Klein surface of genus $p \ge 3$ then
the orientable
two-sheeted covering surface of $Y$ has genus $\ge 2$, so that the
universal covering
 space of $Y$ is the upper half-plane {\cal H} (with the hyperbolic metric)
and the
group of covering transformations is a non-orientable surface subgroup $K$
generated
by glide-reflections. If $G$ is a group of automorphisms of $Y$ then  the
elements of
 $G$ lift to a {\it non-euclidean crystallographic (NEC) group} $\Gamma$
acting on
{\cal H}. There is a smooth epimorphism
\begin{equation}\label{epi}
\theta\co  \Gamma \to G
\end{equation}
 whose kernel is $K$,
where smooth means that $\theta$ preserves the orders of elements of finite
order in
$\Gamma$. The transformation group $(\Gamma, {\cal H})$\, is called the {\it
universal
covering transformation group} of $(G, Y)$.

Now let $G = \< t| t^{2N} = 1\>$ be a cyclic group of order $2N$. As $\theta$ is
 smooth we must have $\theta(c) = t^{N}$ for every reflection $c$ in
$\Gamma$. Also
we cannot have two distinct reflections in $\Gamma$ whose product has
finite order. So
it follows, in the canonical presentation of NEC groups as given in
\cite{[M1]} or  \cite{[BEGG]}, that $\Gamma$ has empty period cycles.

Thus $\Gamma$ has signature of the form
\begin{equation}\label{sign}
s({\Gamma}) = (g; \pm; [m_1,...,m_n]; \{ (\quad)^k \})
\end{equation}
\noindent
with $k$ empty period cycles; then $\Gamma$ has one of the two presentations
 depending on whether there is a $+$ or a $-$ in the signature;

for the $(+)$ case
\begin{eqnarray}\label{pres+}
x_1,\dots , x_n, e_1,\dots ,e_k, c_1,\dots ,c_k, a_1, b_1,\dots ,a_g,b_g \mid
\nonumber \\
x_i^{m_i} = 1, i = 1,...,n,
c_{j}^2 = c_{j}e_j^{-1}c_{j}e_j = 1, j = 1,...,k, \nonumber \\
x_1...x_ne_1...e_ka_1b_1a_1^{-1}b_1^{-1}...a_gb_gagh^{-1}b_g^{-1}
\end{eqnarray}

for the $(-)$ case
\begin{eqnarray}\label{pres-}
x_1,\dots , x_n, e_1,\dots ,e_k, c_1,\dots ,c_k, d_1,...,d_g
 \mid  \nonumber \\  x_i^{m_i} = 1, i = 1,...,n,
c_{j}^2 = c_{j}e_j^{-1}c_{j}e_j = 1, j = 1,...,k,
 x_1...x_ne_1...e_kd_1^2...d_g^2
\end{eqnarray}

In these presentations the generators $x_i$ are elliptic elements, the
generators
$c_j$ are reflections, {\it the generating reflections} of $\Gamma$, and the
 generators $e_j$ are orientation-preserving transformations called the
{\it connecting
generators}. Each empty period cycle corresponds to a conjugacy class of
reflections
in $\Gamma$.

One important fact to note about these presentations is that the connecting
 generator $e_j$ commutes with the generating reflection $c_j$, and in fact
the centralizer
of $c_j$ in $\Gamma$ is just the group $gp\<c_j, \  e_j\> \cong C_2 \times
C_{\infty}$. (See
\cite{[Si2]} )

\section{The fixed-point set of a power of $t$} 

Let $Y$ be a non-orientable surface of topological genus $ p \ge 3$ and let
$t$ be an
automorphism of order $2N$. If $ 1 \le i < 2N$ and $ i\neq
N$ then the number
 of fixed points of the automorphism $t^{i}$ is given by Macbeath's formula
(see
 \cite{[M2]} ). If $t^{i}$ has order $d$ than $t^{i}$ has
\begin{equation}\label{iso}
2N \sum_{d \mid m_j}{1\over m_j}
\end{equation}
fixed points, where $m_j$ runs over the periods in $s(\Gamma )$.

This is because Macbeath's proof (applying to Fuchsian groups) only uses
the facts that
each period corresponds to a unique conjugacy class of elliptic elements of
$\Gamma$,
and each elliptic element has a unique fixed point in {\cal H}. Now, the
number of
 isolated fixed points of $t^{i}$ is independent of the smooth epimorphism
$\theta$ above. However
the epimorphism $\theta$ does play a part in the number of ovals of $t^N$.

\begin{Th}
Let $Y$ be a non-orientable surface of topological genus $ p \ge 3$. Let $G
\cong C_{2N}
 =  \< t \mid t^{2N} = 1 \>$ be a group of automorphisms of $Y$, and let
$\theta$ and
$\Gamma$ be as described in equations \ref{epi} and \ref{sign}. If $\theta
(e_j) =
 t^{v_j}$ than the number of ovals of the involution $t^N$ is
\begin{equation}\label{oval}
\sum_{j=1}^k{ (N, v_j )}
\end{equation}
and the number of isolated fixed points of  $t^N$ is
$$ 2N \sum_{ m_j \  even}{1\over m_j}.$$
\end{Th}

\begin{proof} Let $\Lambda = \theta^{-1} (\<t^N\>)$ so that $\Lambda$
contains the group
$K = Ker\theta$ with index $2$. Now, $\Lambda$ must have signature of the form
\begin{equation}\label{sign2}
s({\Gamma}) = (g; \pm; [2^{(r)}]; \{ (\quad)^s \})
\end{equation}
with $r$ periods equal to $2$ and $s$ empty period cycles.

The reason that all periods in $\Lambda$ are equal to $2$
 is because if $m_j$ in $s(\Gamma)$ is even then
$x_j^{m_j /2} \in \Lambda$ and any elliptic element of $\Lambda$ are
conjugate to
some $x_j^{m_j /2}$ (see \cite{[Si1]} ).

By results in \cite{[BCNS]} (see also \cite{[BEGG]}), r is the number of
isolated
fixed points of $t^N$ and is given by Macbeath's formula
$$ 2N \sum_{ m_j \  even}{1\over m_j}$$

It also follows from \cite{[BCNS]} that the number of ovals of $t^N$ is just
 the number $s$ of period cycles in $\Lambda$, which corresponds to the
number of
 conjugacy classes of reflections in $\Lambda$. As a reflection $c_j$ in
$\Lambda$
belongs also to $\Gamma$ and the group  $\Gamma$ has $k$ conjugacy classes of
 reflections, we just have to determine into how many $\Lambda$--conjugacy
classes the
$\Gamma$--conjugacy class of $c_j$ splits. We shall use the epimorphism
 $\theta$
to calculate this number.

There is a transitive action of $\Gamma$ on the $\Lambda$--conjugacy classes
of $c_j$
in $\Lambda$ by letting $\gamma \in \Gamma$ map the reflection $gc_jg^{-1}$ to
$g\gamma c_j\gamma^{-1}g^{-1}$, with $g\in \Lambda$. (Because $\Lambda
\triangleleft
\Gamma$).
Clearly, if $\lambda \in \Lambda$ then $\lambda$ has a trivial action on these
$\Lambda$--conjugacy classes. So we have an action of $\Gamma / \Lambda
\cong C_{2N} /
C_2 \cong C_N$ on these classes. As the centralizer of $c_j$ in $\Gamma$ is
just
$\< c_j, e_j\>$, the stabilizer of the $\Lambda$--conjugacy classes of $c_j$
in $\Lambda$ are the cosets $\Lambda,\Lambda e_j, \dots, \Lambda
e_j^{\delta_j -1}$,
where $\delta_j = exp_{\Lambda}e_j$, the least positive power of $e_j$ that
belongs to
$\Lambda$. Now, let $\varepsilon_j = exp_K e_j$. Then either $\varepsilon_j
=\delta_j$
or  $\varepsilon_j = 2\delta_j$.

The additive group $Z_{2N}$ contains a subgroup isomorphic to $Z_{N}$ and
$a \in Z_{N}$
has order ${N\over (N,a)}$ in $Z_{N}$ so that $a$ has the same order in
$Z_{2N}$ if and
only if $(2N, a) = 2(N, a)$. If $(2N, a) = (N, a)$ then the order of $a$ in
$Z_{2N}$
is twice the order of $a$ in  $Z_{N}$ and we then find that
$$\varepsilon_j = \delta_j \qquad \mbox{if} \qquad  (2N, v_j) = 2(N, v_j)$$
and
$$\varepsilon_j = 2\delta_j \qquad \mbox{if} \qquad  (2N, v_j) = (N, v_j),$$
where $\theta(e_j) = t^{v_j}$.

By the above argument on the action of $\Gamma / \Lambda$ on the
$\Lambda$--conjugacy
 classes of $c_j$ we see that the number of such classes is $N / \delta_j$,
which is

\noindent
if $\varepsilon_j = \delta_j$
$${N \over \delta_j}\quad = \quad {N \over \varepsilon_j}\quad = \quad
{N(2N, v_j)
\over2N } \quad = \quad {(2N, v_j)\over2 }\quad = \quad (N, v_j),$$

\noindent
or if $\varepsilon_j = 2\delta_j$
$${N \over \delta_j}\quad = \quad {2N \over \varepsilon_j}\quad = \quad
{2N(2N, v_j)
\over2N } \quad = \quad (2N, v_j) \quad = \quad (N, v_j)$$
Thus in both cases the generating reflection $c_j$ of $\Gamma$ induces $(N,
v_j)$ conjugacy
classes of reflections in $\Lambda$. Thus the number of ovals of $t^{N}$ in
$Y$ is
\begin{equation}
\sum_{j=1}^k{ (N, v_j )}
\end{equation}
\end{proof}

\begin{Th}
The ovals of $t^{N}$ in $Y$ induced by the $j$th period cycle in $\Gamma$
are twisted
if $(2N, v_j) = (N, v_j)$ and untwisted if $(2N, v_j) = 2(N, v_j)$.
\end{Th}

\begin{proof}As we have found in Theorem 3.1, the $j$th empty period cycle in
$\Gamma$ induces $(N,   
v_j)$ empty period
cycles in $\Lambda$.  The generating
reflections of these
period cycles are just conjugates of $c_j$ in $\Gamma$ and, as the
corresponding
 connecting generator $e_j$ is just the orientation-preserving element
generating the
centralizer of $c_j$ in $\Gamma$, we see that the connecting generator of
each of the
period cycles in $\Lambda$ induced by the $j$th period cycle in $\Gamma$ is
just
conjugate to $e_j^{\delta_j}$, $\delta_j = exp_{\Lambda}e_j$ as in the
proof of Theorem 3.1.
Now, let $\theta' \co  \Lambda \to C_2 = gp\<\xi \>$, where $\xi = t^N$, be the
restriction
of the epimorphism $\theta \co  \Gamma \to C_{2N}$. Then

\noindent
if $\varepsilon_j = \delta_j$
$$\theta'(e_j^{\delta_j})\quad = \quad \theta'(e_j^{\varepsilon_j}) \quad =
\quad
\theta(e_j^{\varepsilon_j}) \quad = \quad 1$$

\noindent
if $\varepsilon_j = 2\delta_j$
$$\theta'(e_j^{\delta_j})\quad = \quad \theta'(e_j^{\varepsilon_j \over 2})
\quad = \quad
\theta(e_j^{\varepsilon_j \over 2}) \quad = \quad \xi ,$$
$\xi$ the generator of $C_2$. Generally, if $c$ is the generating
reflection of an empty
 period cycle of $\Lambda$ and $e$ is the corresponding connecting
generator then figures 1 and 2 show 
that $\theta'(e) =1$ corresponds to an untwisted
oval while
$\theta'(e) = \xi$ corresponds to a twisted oval.

\begin{figure}[htb]       
\cl{\relabelbox\small\epsfxsize 2.5truein\epsfbox{figor.ai}
\relabel{ep}{$\varepsilon$}
\relabel{epp}{$\varepsilon'$}
\relabel{cep}{$c(\varepsilon)$}
\relabel{cepp}{$c(\varepsilon')$}
\relabel{F}{$F$}
\relabel{cF}{$c(F)$}
\relabel{g}{$\gamma$}
\endrelabelbox
\relabelbox\small\epsfxsize 2.5truein\epsfbox{fignonor.ai}
\relabel{ep}{$\varepsilon$}
\relabel{epp}{$\varepsilon'$}
\relabel{cep}{$c(\varepsilon)$}
\relabel{cepp}{$c(\varepsilon')$}
\relabel{F}{$F$}
\relabel{cF}{$c(F)$}
\relabel{g}{$\gamma$}
\endrelabelbox}
\caption{$\theta'(e) =1$ so $e \in K$\qquad\qquad Figure 2:\stdspace
$\theta'(e) = \xi$ so $ce \in K$}
\end{figure}

However, as in the proof of Theorem 3.1 $\varepsilon_j = \delta_j$ if and
only if
$(2N, v_j) = 2(N, v_j)$ and hence we have untwisted ovals while
$\varepsilon_j = 2\delta_j$ if and only if
$(2N, v_j) = (N, v_j)$ and we have twisted ovals. \end{proof}  

\section{Bounds and examples}

In \cite{[Sch]} (also see \cite{[BCNS]})  Scherrer showed that that if an
involution   of a
non-orientable surface of  genus $p$ has $\mid F \mid$ fixed points and $\mid V
\mid$ ovals then
$$\mid F \mid +2\mid V \mid\le p+2.$$ In our examples we will show that for
any integer
$N$ we can find a non-orientable surface of genus $p$ admitting a $C_{2N}$
action with generator $t$
such that $t^N$ attains the Scherrer bound.

{\bf Example 1}\stdspace Bujalance \cite{[B]} found the maximum order for an
automorphism $t$ of
a non-orientable surface $Y$ of genus $p \ge 3$; it is $2p$ for odd $p$ and
$2(p-1)$
 for even $p$. The universal covering transformation group $\Gamma$ has
signature
$ s(\Gamma) = (0; [2, p]; \{ (\quad)\})$ for odd $p$, and signature
$s(\Gamma) =
 (0; [2,  2(p-1)]; \{ (\quad)\})$ for even $p$. There is, essentially, only
one way of
defining the epimorphism $\theta$ in each case:

if $p$ is odd, we define  $\theta \co  \Gamma \to C_{2p}$ by $\theta(x_1) = t^p$,
 $\theta(x_2) = t^2$,  $\theta(c) = t^p$, and  $\theta(e) = t^{p-2}$,

if $p$ is even, we define  $\theta \co  \Gamma \to C_{2(p-1)}$ by $\theta(x_1)
= t^{p-1}$,
 $\theta(x_2) = t^1$,  $\theta(c) = t^{p-1}$, and  $\theta(e) = t^{p-2}$.

Using Macbeath's formula (\ref{iso}) we see that the involution $t^p$ has $p$
 fixed points for surfaces of both odd and even genera.
 Now, if $p$ is odd then the involution $t^p$ also has,
by Theorems
3.1 and 3.2, one twisted oval if $p$ is odd as
$(p, p-2) = (2p, p-2) = 1$. If $p$ is even then the involution $t^{p-1}$
has, by
Theorems 3.1 and 3.2, one untwisted oval as
$(p-1, p-2) = 1$ and $(2(p-1), p-2) = 2(p, p-2) = 2$. We note that the
involution $t^p$
obeys the Scherrer bound. Note that the orders of the cyclic groups in
Bujulance's examples are $\equiv 2$ mod 4. Our
second example shows that the Scherrer bound can be obtained for the
involution in a $C_4$ action.

{\bf Example 2}\stdspace Let $Y$ be a non-orientable surface of  genus $p\ge 3$,
and let
$t$ be an automorphism of $Y$ of order $4$. Let $\Gamma$ have signature
$$(0;+;[2^{(r)},4,4]; (\quad)^k)$$
and define a smooth epimorphism $\theta\co \Gamma\rightarrow C_4$ by mapping
the generators of order
two to $t^2$, the two generators of order 4 to $t$ and $t^{-1}$ and the
connecting generators to
the identity. We then find that for the involution $t^2$,  $\mid F
\mid=2r+2$, and $\mid V
\mid=2k$,and $p=4k+2r$, so that we find infinitely many surfaces where the Scherrer
bound is attained for the
involution in $C_4$. This is easily extended to groups of order $4m$ by
replacing the two periods 4
in the signature of $\Gamma$ by $4m$.

\Addresses\recd

\end{document}